\documentclass[12pt]{amsart}
\usepackage{amssymb}
\makeatletter
\renewcommand{\section}{\@startsection{section}{1}%
  \z@{.7\linespacing\@plus\linespacing}%
  {-.5em}%                     % si < 0 pas de passage a la ligne
  {\bfseries\large}}           % Le style (taille, bold ...)

\renewcommand{\subsection}{\@startsection{subsection}{2}%
  \z@{.5\linespacing\@plus.7\linespacing}{-.5em}%
  {\normalfont\bfseries}}

\renewcommand{\subsubsection}{\@startsection{subsubsection}{3}%
  \z@{.5\linespacing\@plus.7\linespacing}{-.5em}%
  {\normalfont\itshape}}
\makeatother

\newenvironment{Prf}{\vskip 1em{\it Proof} :}%
{\unskip\hfill\null\nobreak\hfill\carre\vskip1em\par}
\newcommand{\carre}{\rule{1ex}{1ex}}

\def \rb {\Bbb{R}}
\def \nb {\Bbb{N}}

\begin{document}

{\centerline{\bf Infinitely many homoclinic orbits for a class of}}
{\centerline{\bf superquadratic Hamiltonian systems}}
$$   $$
\centerline{ Mohsen TIMOUMI. Dpt of Mathematics}
\centerline{ Faculty of Sciences. 5019 Monastir. Tunisia}
\centerline{email:m\_timoumi@yahoo.com}
$$   $$
{\bf Abstract.} In this paper, we prove the existence of infinitely many homoclinic orbits for the first order Hamiltonian systems $J\dot{x}-M(t)x+ R'(t,x)=0$, by the minimax methods in critical point theory, when $R(t,y)$ satisfies the superquadratic condition ${{R(t,x)}\over{\left|x\right|^{2}}}\longrightarrow \pm\infty$ as $\left|x\right|\longrightarrow\infty$, uniformly in $t$, and need not satisfy the global Ambrosetti-Rabinowitz condition.
\section{Introduction and main result.} In this paper, we are interested in the existence of homoclinic solutions for the first-order Hamiltonian system
$$\dot{x}=JH'(t,x) \leqno(\mathcal H)$$
where $x(t)=(p(t),q(t))\in\rb^{N}\times\rb^{N}$, $J=\left(\begin{array}{ll}
0 & -I_{N}\\
I_{N} & 0
\end{array}\right)$ is the standard symplectic matrix and $H$ is of the type
$$H(t,x)=-\frac{1}{2}M(t)x.x +R(t,x)$$
with $M(t)=\left(\begin{array}{ll}
0 & L(t)\\
L(t) & 0
\end{array}\right)$, $L$ is a continuous function taking values in the set of $(N\times N)-$symmetric matrices and $R:\rb\times\rb^{2N}\longrightarrow\rb$, $(t,x)\longmapsto R(t,x)$ is a continuous function, differentiable with respect to the second variable with continuous derivative $R'(t,x)=\frac{\partial R}{\partial x}(t,x)$. Here $x.y$ denotes the Euclid's inner product of $x,y\in\rb^{2N}$ and $\left|.\right|$ denotes the corresponding Euclid's norm. As usual, assuming that $x=0$ is an equilibrium for $(\mathcal H)$, we say that a solution $x$ of $(\mathcal H)$ is homoclinic to $0$ if $x\in C^{1}(\rb,\rb^{2N})$ satisfies $x\neq 0$ and the asymptotic condition $x(t)\longrightarrow 0$ as $\left|t\right|\longrightarrow\infty$.\\
Establishing the existence of homoclinic orbits of Hamiltonian systems is one of the most important problem in the theory of Hamiltonian systems. During the two last decades, the existence and multiplicity of homoclinic solutions for Hamiltonian systems have been extensively investigated by many authors with the aid of the variational methods. For examples see [1-10] for the second-order systems, and  [11-17] for the first-order systems. \\
All classical known results for Hamiltonian systems (see [1-3,8,11-16]) are obtained under the following assumption that the Hamiltonian satisfies the so-called Ambrosetti-Rabinowitz condition, that is, there exists a constant $\mu>2$ such that for all $t\in \rb$ and $x\in\rb^{2N}$, $\left|x\right|\geq r$,
$$0<\mu R(t,x)\leq R'(t,x).x.$$
In very recent years, many authors devoted to the existence of homoclinic orbits for second order systems (see [5-7, 9,10]) under a kind of new superquadratic conditions firstly introduced by Fei [19] for the existence of periodic solutions. Motivated by the work of [20], the author gets recently in [17], without the Ambrosetti-Rabinowitz condition, the existence of at least one homoclinic orbit for $(\mathcal H)$ under the following conditions:\\
$(L_{1})$ There exists a constant $\gamma<0$ such that the smallest eigenvalue
$$l(t)=\inf_{\left|\xi\right|=1}L(t)\xi.\xi$$
of $L(t)$ satisfies
$$l(t)\left|t\right|^{\gamma-1}\longrightarrow\infty\ as\ \left|t\right|\longrightarrow\infty;$$
$(L_{2})$ $L\in C^{1}(\rb,\rb^{2N})$ and there is $T_{0}>0$ such that $2L(t)+\dot{L}(t)$ and 
$2L(t)-\dot{L}(t)$ are nonnegative definite for all $\left|t\right|\geq T_{0}$, where $\dot{L}(t)=\frac{dL}{dt}(t)$; \\
$${{R(t,x)}\over{\left|x\right|^{2}}}\longrightarrow+\infty\ as\ \left|x\right|\longrightarrow\infty,\ uniformly\ in\ t\in\rb, \leqno(R_{1})$$
$${{\left|R'(t,x)\right|}\over{\left|x\right|}}\longrightarrow 0\ as\ \left|x\right|\longrightarrow 0,\ uniformly\ in\ t\in\rb; \leqno(R_{2})$$
$(R_{3})$ there exist $a>0$ and $\alpha>1$ such that
$$\left|R'(t,x)\right|\leq a(\left|x\right|^{\alpha}+1),\ \forall t\in\rb,\ \forall x\in\rb^{2N};$$
$(R_{4})$ there exist $\beta>\alpha$, $b>0$ and $r>0$ such that
$$R'(t,x).x-2R(t,x)\geq b\left|x\right|^{\beta},\ \forall t\in\rb,\ \forall \left|x\right|\geq r,$$
$$R'(t,x).x\geq 2R(t,x)\geq 0,\ \forall t\in\rb,\ \forall x\in\rb^{2N}. \leqno(R_{5})$$
{\bf Remark 1.1.} It is easy to see that the function
$$R(t,x)=\left|sin(t)\right|\left|x\right|^{\frac{3}{2}}+\left|x\right|^{2}ln(1+\left|x\right|^{2}),\ \forall t\in\rb,\ \forall x\in\rb^{2N}$$
satisfies $(R_{1})-(R_{5})$ and don't satisfies the Ambrosetti-Rabinowitz condition.\\
The purpose of this paper is to show that the Hamiltonian system $(\mathcal{H})$ possesses infinitely many homoclinic orbits if $R(t,x)$ is even in $x$ and satisfies the above assumptions. \\
Our main result reads as follows:\\
{\bf Theorem 1.1.} Assume $(L_{1})$, $(L_{2})$ and $(R_{1})-(R_{5})$ hold, and suppose, in addition, that $R(t,-x)=R(t,x)$ for all $(t,x)\in\rb\times\rb^{2N}$. Then $(\mathcal{H})$ possesses infinitely many homoclinic orbits $(x_{k})$ such that
$$\int_{\rb}[-{1\over 2}J\dot{x}_{k}.x_{k}+{1\over 2}M(t)x_{k}.x_{k}-R(t,x_{k})]dt\longrightarrow +\infty\ as\ k\longrightarrow\infty.$$
{\bf Remark 1.2.} Observe that if $x$ is a solution of $(\mathcal{H})$ then $y(t)=x(-t)$ is a solution of the system
$$J\dot{y}(t)+M(-t)y(t)-R'(-t,y(t))=0.$$
Moreover, $-R(-t,x)$ satisfies $(R_{1})$, $(R_{4})$, $(R_{5})$ whenever $R(t,x)$ satisfies respectively the following assumptions
$${{R(t,x)}\over{\left|x\right|^{2}}}\longrightarrow-\infty\ as\ \left|x\right|\longrightarrow\infty,\ uniformly\ in\ t\in\rb; \leqno(R'_{1})$$
$(R'_{4})$ there exist $\beta>\alpha$, $b>0$ and $r>0$ such that
$$R'(t,x).x-2R(t,x)\leq -b\left|x\right|^{\beta},\ \forall t\in\rb,\ \forall \left|x\right|\geq r;$$
$$R'(t,x).x\leq 2R(t,x)\leq 0,\ \forall t\in\rb,\ \forall x\in\rb^{2N}. \leqno(R'_{5})$$
So we have\\
{\bf Theorem 1.2.} Assume $(L_{1})$, $(L_{2})$ and $(R'_{1})$, $(R_{2})$, $(R_{3})$, $(R'_{4})$, $(R'_{5})$ hold, and suppose, in addition, that $R(t,-x)=R(t,x)$ for all $(t,x)\in\rb\times\rb^{2N}$. Then $(\mathcal{H})$ possesses infinitely many homoclinic orbits $(x_{k})$ such that
$$\int_{\rb}[-{1\over 2}J\dot{x}_{k}.x_{k}+{1\over 2}M(t)x_{k}.x_{k}-R(t,x_{k})]dt\longrightarrow -\infty\ as\ k\longrightarrow\infty.$$
\section{Preliminary results.} For a selfadjoint operator $S$ acting in $L^{2}=L^{2}(\rb,\rb^{2N})$, let $\mathcal{D}(S)$ denotes the domain, $\left|S\right|$ the absolute value, and $\left|S\right|^{\frac{1}{2}}$ the square root. In the following, $c_{i}$ denotes a positive constant. If $S_{1}$ and $S_{2}$ are two selfadjoint operators with $\mathcal{D}(S_{1})\subset \mathcal{D}(S_{2})$, then $\mathcal{D}(\left|S_{1}\right|^{\frac{1}{2}})\subset \mathcal{D}(\left|S_{2}\right|^{\frac{1}{2}})$, and if moreover $\left\|S_{1}u\right\|\leq c_{1}\left\|S_{2}u\right\|$ for all $u\in\mathcal{D}(S_{1})$, then $\left\|\left|S_{1}\right|^{\frac{1}{2}}u\right\|_{L^{2}}\leq c_{2}\left\|\left|S_{2}\right|^{\frac{1}{2}}u\right\|_{L^{2}}$ for all $u\in \mathcal{D}(\left|S_{1}\right|^{\frac{1}{2}})$.\\
Let $J_{0}$ be the $(2N)\times(2N)$ matrix
$$J_{0}=\left(\begin{array}{ll}
0 & I_{N}\\
I_{N} & 0
\end{array}\right),$$
and $A_{0}=-J\frac{d}{dt}+J_{0}$, a selfadjoint operator acting in $L^{2}$. Let $W^{1,2}=W^{1,2}(\rb,\rb^{2N})$ be the usual Sobolev space. Clearly, $<A_{0}u,A_{0}u>_{L^{2}}=\left\|u\right\|^{2}_{W^{1,2}}$ for all $u\in\mathcal{D}(A_{0})$, and $\mathcal{D}(A_{0})=W^{1,2}$, where $<.,.>_{L^{2}}$ denotes the $L^{2}$ inner product.  Note that letting $D=(-I_{2N}\frac{d^{2}}{dt^{2}}+I_{2	N})^{\frac{1}{2}}$, one has $\mathcal{D}(D)=W^{1,2}$ and $\mathcal{D}(\left|D\right|^{\frac{1}{2}})=W^{\frac{1}{2},2}(\rb,\rb^{2N})=H^{\frac{1}{2}}$, the Sobolev space of fractional order. Hence $\mathcal{D}(\left|A_{0}\right|^{\frac{1}{2}})=H^{\frac{1}{2}}$ and
$$\left\|u\right\|_{H^{\frac{1}{2}}}\leq c_{3}\left\|\left|A_{0}\right|^{\frac{1}{2}}u\right\|_{L^{2}},\ \forall u\in\mathcal{D}(\left|A_{0}\right|^{\frac{1}{2}}). \leqno(2.1)$$
Suppose that $L$ satisfies $(L_{1})$ and $(L_{2})$. Let $A=-J\frac{d}{dt}+ M$, a selfadjoint operator with $\mathcal{D}(A)\subset L^{2}$, defined as a sum of quadratic forms. Let $\left\{E(\lambda)/-\infty<\lambda<\infty\right\}$ denotes the resolution of $A$, and $U=I-E(0)-E(-0)$.  Then $U$ commutes with $A$, $\left|A\right|$ and $\left|A\right|^{1\over 2}$, and $A=\left|A\right|U$ is the polar decomposition of $A$ (see [18]). $\mathcal{D}(A)$ is a Hilbert space equipped with the norm
$$\left\|u\right\|_{1}=\left\|(I+\left|A\right|)u\right\|_{L^{2}},\ \forall u\in\mathcal{D}(A).$$
It is easy to verify 
$$\left\|u\right\|_{W^{1,2}}\leq c_{4}\left\|u\right\|_{1},\ \forall u\in \mathcal{D}(A),\leqno(2.2)$$
i.e., $\mathcal{D}(A)$ is continuously embedded in $W^{1,2}$. Combining (2.1) and (2.2) yields
$$\left\|u\right\|_{H^{\frac{1}{2}}}\leq c_{5}\left\|(I+\left|A\right|)^{\frac{1}{2}}u\right\|_{L^{2}},\ \forall u\in\mathcal{D}((I+\left|A\right|)^{\frac{1}{2}}). \leqno(2.3)$$
Moreover $\mathcal{D}(A)$ is compactly embedded in $L^{2}(\rb,\rb^{2N})$ (see [14]). Then $(I+\left|A\right|)^{-1}: L^{2}\longrightarrow L^{2}$ is a compact linear operator. Therefore a standard argument shows that the spectrum $\sigma(A)$ of $A$ consists of eigenvalues numbered by (counted in their multiplicities)
$$...\lambda_{-2}\leq\lambda_{-1}\leq 0<\lambda_{1}\leq\lambda_{2}...$$
with $\lambda_{k}\longrightarrow \pm\infty$ as $k\longrightarrow\pm\infty$, and a corresponding system of eigenfunctions $(e_{k})$ of $A$ forms an orthonormal basis in $L^{2}$. \par\noindent 
Now set $E=\mathcal{D}(\left|A\right|^{1\over 2})=\mathcal{D}((I+\left|A\right|)^{1\over 2})$. $E$ is a Hilbert space under the inner product
$$<u,v>_{0}=<\left|A\right|^{1\over 2}u,\left|A\right|^{1\over 2}v>_{L^{2}}+<u,v>_{L^{2}}$$
and norm
$$\left\|u\right\|_{0}=<u,u>^{1\over 2}_{0}=\left\|(I+\left|A\right|)^{1\over 2}u\right\|_{L^{2}}.$$
Let $E^{0}=Ker(A)$, $E^{+}=CL_{E}(span\left\{e_{1},...,e_{n}\right\})$ and $E^{-}=(E^{0}\oplus E^{+})^{\bot_{E}}$, where $CL_{E}S$ stands for the closure of $S$ in $E$ and $S^{\bot_{E}}$ the orthogonal complementary subspace of $S$ in $E$. Then
$$E=E^{-}\oplus E^{0}\oplus E^{+}$$
is an orthonormal decomposition of $E$. Since $\mathcal{D}(A)$ is compactly embedded in $L^{2}$, $0$ is at most an isolated eigenvalue of $A$, then for the later convenience, we introduce on $E$ the following inner product
$$<u,v>=<\left|A\right|^{1\over 2}u,\left|A\right|^{1\over 2}v>_{L^{2}}+<u^{0},v^{0}>_{L^{2}}$$
for all $u=u^{-}+u^{0}+u^{+},\ v=v^{-}+v^{0}+v^{+} \in E^{-}\oplus E^{0}\oplus E^{+}$, and norm
$$\left\|u\right\|=<u,u>^{\frac{1}{2}}.$$
Since $\gamma<0$, then $\frac{2}{2-\gamma}<1$ and we get\\
{\bf Proposition 2.1.}[14] Let $L$ satisfies $(L_{1})$ and $(L_{2})$. Then $E$ is compactly embedded in $L^{p}$ for all $[1,\infty[$, which implies that for all $p\in[1,\infty[$, there exists a constant $\lambda_{p}>0$ such that
$$\left\|u\right\|_{L^{p}}\leq \lambda_{p}\left\|u\right\|,\ \forall u\in E. \leqno(2.1)$$
Finally, we introduce
$$a(u,v)=<\left|A\right|^{1\over 2}U u,\left|A\right|^{1\over 2}v>_{L^{2}} \leqno(2.2)$$
for all $u,v\in E$. $a(u,u)$ is the quadratic form associated with $A$. Clearly, for $u\in\mathcal{D}(A)$ and $v\in E$ we have
$$a(u,v)=<Au,v>_{L^{2}}=\int_{\rb}(-J\dot{u}+M u).u dt. \leqno(2.3)$$
Plainly, $E^{-}$, $E^{0}$ and $E^{+}$ are orthogonal to each other with respect to $a$, and moreover
$$a(u,v)=<(P^{+}-P^{-})u,v>,\ \forall u,v\in E,$$
$$a(u,u)=\left\|u^{+}\right\|^{2}-\left\|u^{-}\right\|^{2},\ \forall u\in E, \leqno(2.4)$$ 
where $P^{\pm}:E\longrightarrow E^{\pm}$ are the orthogonal projectors and $u=u^{-}+u^{0}+u^{+}\in E^{-}\oplus E^{0}\oplus E^{+}$.\\
The following critical point proposition will be used for proving the previous Theorem. \\
Let $E$ be a real Hilbert space with the norm $\left\|.\right\|$. Suppose that $E$ has an orthogonal decomposition $E=E^{1}\oplus E^{2}$ with both $E^{1}$ and $E^{2}$ being infinite dimensional. Suppose $(v_{n})$ (resp. $(w_{n})$) is an orthonormal basis for $E^{1}$ (resp. $E^{2}$), and set
$$X_{n}=span\left\{v_{1},...,v_{n}\right\}\oplus E^{2},\ X^{m}=E^{1}\oplus\left\{w_{1},...,w_{m}\right\}.$$
Recall that we say $f\in C^{1}(E,\rb)$ satisfies $(PS)^{*}$ condition if any sequence $(x_{n})$ with $x_{n}\in X_{n}$ for which $0\leq f(x_{n})\leq const.$ and $f^{'}_{n}(x_{n})\longrightarrow 0$ as $n\longrightarrow\infty$ possesses a convergent subsequence, where $f_{n}=f_{/X_{n}}$. We also say that $f$ satisfies $(PS)^{**}$ condition if for each $n\in\nb$, $f_{n}$ satisfies the Palais-Smale condition, i.e., any sequence $(u_{_{k}})\subset X_{n}$ for which $(f(u_{k}))$ is bounded and $f^{'}_{n}(u_{k})\longrightarrow 0$ as $k\longrightarrow\infty$ possesses a convergent subsequence. \\
{\bf Proposition 2.2.}[11] Let $E$ be as above and let $f\in C^{1}(E,\rb)$ be even, satisfy $(PS)^{*}$ and $(PS)^{**}$, and $f(0)=0$. Suppose moreover that $f$ satisfies, for any $m\in\nb$, \par\noindent
$(f_{1})$ there is $R_{m}>0$ such that $f(u)\leq 0,\ \forall u\in X^{m}$ with $\left\|u\right\|\geq R_{m}$; \\
$(f_{2})$ there are $r_{m}>0$, $a_{m}>0$ with $a_{m}\longrightarrow\infty$ as $m\longrightarrow\infty$ such that
$$f(u)\geq a_{m},\ \forall u\in (X^{m-1})^{\bot}\ with\ \left\|u\right\|=r_{m};$$
$(f_{3})$ $f$ is bounded from above on bounded sets of $X^{m}$. \\
Then $f$ has a positive critical value sequence $(c_{k})$ satisfying $c_{k}\longrightarrow\infty$ as $k\longrightarrow\infty$. \\
\section{Proof of Theorem 1.1.} 
Define a functional $f$ in $E$ by
$$f(u)={1\over 2}a(u,u)-\int_{\rb}R(t,u)dt,\ \forall u\in E. \leqno(3.1)$$
By (2.4), we have
$$f(u)={1\over 2}(\left\|u^{+}\right\|^{2}-\left\|u^{-}\right\|^{2})-\int_{\rb}R(t,u)dt, \leqno(3.2)$$
for all $u=u^{-}+u^{0}+u^{+}\in E^{0}\oplus E^{-}\oplus E^{+}$. \\
By $(R_{2})$, $(R_{3})$, for all $\epsilon >0$ there exists a constant $C_{\epsilon}>0$ such that
$$\left|R'(t,x)\right|\leq 2\epsilon\left|x\right|+(\alpha+1)C_{\epsilon}\left|x\right|^{\alpha},\ \forall t\in\rb,\ \forall x\in \rb^{2N}. \leqno(3.3)$$
Now, using the Mean Value Theorem and the inequality and (3.3) imply
$$\left|R(t,x)\right|=\left|\int^{1}_{0} R'(t,sx).x ds\right|$$
$$\leq\epsilon\left|x\right|^{2}+C_{\epsilon}\left|x\right|^{\alpha+1},\ \forall t\in\rb,\ \forall x\in\rb^{2N}. \leqno(3.4)$$
It is well known that (3.3) and (3.4) imply that the functional $f$ is continuously differentiable in $E$ and for all $u,v\in E$
$$f'(u)v=\int_{\rb}(-J\dot{u}+M(t)u-R'(t,u)).v dt. \leqno(3.5)$$
Moreover, the critical points of $f$ on $E$ are exactly the homoclinic orbits of the system $(\mathcal H)$.\\
Now, let $E^{1}=E^{-}\oplus E^{0}$ and $E^{2}=E^{+}$ with $(v_{n}=e_{n})^{\infty}_{1}$ and $(w_{m}=e_{-m})^{\infty}_{1}$ respectively, where $(e_{n})^{\infty}_{-\infty}$ is the system of eigenfunctions of $A$ (see section 2). Set also $X_{n}=span\left\{v_{1},...,v_{n}\right\}\oplus E^{2}$, $X^{m}=E^{1}\oplus\left\{w_{1},...,w_{m}\right\}$ and $f_{n}=f_{/X_{n}}$.
We will verify that $f$ satisfies the assumptions of Proposition 2.2. We will proceed by successive lemmas.\\
{\bf Lemma 3.1.} $f$ satisfies $(f_{1})$.\\
\begin{Prf} We claim that there exists $\epsilon_{1}>0$ such that
$$meas(\left\{t\in\rb/\left|u(t)\right|\geq\epsilon_{1}\left\|u\right\|\right\})\geq\epsilon_{1},\ \forall u\in X^{m}-\left\{0\right\}. \leqno(3.6)$$
Otherwise, for any positive integer $k$, there exists $u_{k}\in X^{m}-\left\{0\right\}$ such that
$$meas(\left\{t\in\rb/\left|u_{k}(t)\right|\geq{1\over k}\left\|u_{k}\right\|\right\})<{1\over k}$$
for all positive integer $k$. Replacing $u_{k}$ by ${u_{k}}\over{\left\|u_{k}\right\|}$, if necessary, we may assume that $\left\|u_{k}\right\|=1$ and
$$meas(\left\{t\in\rb/\left|u_{k}(t)\right|\geq {1\over k}\right\})<{1\over k} \leqno(3.7)$$
for all positive integer $k$. Since $X^{m}$ is closed in $H^{1\over 2}$, it is reflexive and then there exists a subsequence, denoted by $(u_{k})$, such that $(u_{k})$ converges weakly to some $u_{0}$ in $X^{m}$. Hence by Proposition 2.1, we can assume, by going to a subsequence if necessary, that $u_{k}\longrightarrow u_{0}$ in $L^{2}$, i.e.
$$\int_{\rb}\left|u_{k}-u_{0}\right|^{2}dt\longrightarrow 0\ as\ k\longrightarrow\infty. \leqno(3.8)$$
Thus there exist $\delta_{1}>0$, $\delta_{2}>0$ such that
$$meas(\left\{t\in\rb/\left|x_{0}(t)\right|\geq \delta_{1}\right\})\geq\delta_{2}. \leqno(3.9)$$
In fact, if not, we have
$$meas(\left\{t\in\rb/\left|u_{0}(t)\right|\geq {1\over k}\right\})=0,$$
for all positive integer $k$, which implies that
$$0\leq\int_{\rb}\left|u_{0}\right|^{4}dt\leq\left\|u_{0}\right\|^{2}_{L^{\infty}}\int_{\rb}\left|u_{0}\right|^{2}dt\leq{{\lambda^{2}_{2}}\over{k^{2}}}\left\|u_{0}\right\|^{2}={{\lambda^{2}_{2}}\over{k^{2}}}\longrightarrow 0$$
as $k\longrightarrow\infty$ by Proposition 2.1. Hence $u_{0}=0$, which contradicts that $\left\|u_{0}\right\|=1$. Therefore (3.9) holds. \\
Now let
$$I_{0}=\left\{t\in\rb/\left|u_{0}(t)\right|\geq \delta_{1}\right\},\ I_{k}=\left\{t\in\rb/\left|u_{k}(t)\right|< {1\over k}\right\},$$
and
$$I^{c}_{k}=\rb-I_{k}.$$
By (3.7) and (3.9), we have for all positive integer $k$
$$meas(I_{k}\cap I_{0})=meas(I_{0}-(I^{c}_{k}\cap I_{0}))\geq meas(I_{0})-meas(I^{c}_{k}\cap I_{0})\geq \delta_{2}-{1\over k}.$$
Let $k$ be large enough such that
$$\delta_{1}-{1\over k}\geq{1\over 2}\delta_{1}\ and\ \delta_{2}-{1\over k}\geq{1\over 2}\delta_{2},$$
one has
$$\left|u_{k}(t)-u_{0}(t)\right|^{2}\geq \left|\left|u_{k}(t)\right|-\left|u_{0}(t)\right|\right|^{2}\geq(\delta_{1}-{1\over k})^{2}\geq{1\over 4}\delta^{2}_{1},\ \forall t\in I_{k}\cap I_{0}$$
which implies that
$$\int_{\rb}\left|u_{k}(t)-u_{0}(t)\right|^{2}dt\geq\int_{I_{k}\cap I_{0}}\left|u_{k}(t)-u_{0}(t)\right|^{2}dt\geq {1\over 4}\delta^{2}_{1}meas(I_{k}\cap I_{0})$$
$$\geq{1\over 4}\delta^{2}_{1}(\delta_{2}-{1\over k})\geq{1\over 8}\delta^{2}_{1}\delta_{2}>0$$
for all large integer $k$. This is a contradiction with (3.8). Therefore (3.6) holds. \\
For $u=u^{-}+u^{0}+u^{+}\in X^{m}$, let
$$\Omega_{u}=\left\{t\in\rb/\left|u(t)\right|\geq\epsilon_{1}\left\|u\right\|\right\}.$$
By $(R_{1})$, for $d={1\over{2\epsilon^{3}_{1}}}>0$, there exists $r_{1}>0$ such that
$$R(t,x)\geq d\left|x\right|^{2},\ \forall \left|x\right|\geq r_{1},\ \forall t\in\rb.$$
Hence one has
$$R(t,u(t))\geq d\left|u(t)\right|^{2}\geq{1\over{2\epsilon_{1}}}\left\|u\right\|^{2} \leqno(3.10)$$
for all $u\in X^{m}$ with $\left\|u\right\|\geq{{r_{1}}\over{\epsilon_{1}}}$ and $t\in \Omega_{u}$. It follows from $(R_{5})$, (3.6) and (3.10) that
$$f(u)={1\over 2}\left\|u^{+}\right\|^{2}-{1\over 2}\left\|u^{-}\right\|^{2}-\int_{\rb}R(t,u)dt$$
$$\leq{1\over 2}\left\|u^{+}\right\|^{2}-\int_{\Omega_{u}}R(t,u)dt
\leq{1\over 2}\left\|u^{+}\right\|^{2}-{1\over{2\epsilon_{1}}}\left\|u\right\|^{2}meas(\Omega_{u})$$
$$\leq{1\over 2}\left\|u^{+}\right\|^{2}-{1\over 2}\left\|u\right\|^{2}\leq 0, \leqno(3.11)$$
for all $u\in X^{m}$ with $\left\|u\right\|\geq \frac{R_{1}}{\epsilon_{1}}$. The proof of Lemma 3.1 is complete.\\
\end{Prf}
{\bf Lemma 3.2.} $f$ satisfies $(f_{2})$.\\
\begin{Prf} Define
$$\eta_{m}=\sup_{u\in(X^{m})^{\bot}-\left\{0\right\}}{{\left\|u\right\|_{L^{\alpha+1}}}\over{\left\|u\right\|}}.$$
Clearly, $\eta_{m}\geq\eta_{m+1}>0$. We claim that
$$\eta_{m}\longrightarrow 0\ as\ m\longrightarrow\infty. \leqno(3.12)$$
Arguing indirectly, assume $\eta_{m}\longrightarrow \eta>0\ as\ m\longrightarrow\infty$. Then there is a sequence $u_{m}\in(X)^{\bot}_{m}$ with $\left\|u_{m}\right\|=1$ and $\left\|u_{m}\right\|_{L^{\alpha+1}}\geq{\eta\over 2}$. For each $w_{k}$, there is $m_{0}\in\nb$ such that
$$\forall m\geq m_{0},\ w_{k}\in X^{m},\ and\ then\ <u_{m},w_{k}>=0.$$
Therefore $<u_{m},w_{k}>\longrightarrow 0$ as $m\longrightarrow\infty$ and then $u_{m}\rightharpoonup 0$ weakly in $E$. So, by Proposition 2.1, $\left\|u_{m}\right\|_{L^{\alpha+1}}\longrightarrow 0$ as $m\longrightarrow\infty$, contradiction. Thus $\eta_{m}\longrightarrow 0$ as $m\longrightarrow\infty$.\par\noindent
By (3.4) with $\epsilon={1\over{4\lambda^{2}_{2}}}$ and $c=C_{\epsilon}$ one has, for $u\in(X^{m-1})^{\bot}$
$$f(u)={1\over 2}\left\|u\right\|^{2}-\int_{\rb}R(t,u)dt\geq{1\over 2}\left\|u\right\|^{2}-\int_{\rb}[\frac{1}{4\lambda^{2}_{2}}\left|u\right|^{2}+c\left|u\right|^{\alpha+1}]dt$$
$$\geq{1\over 2}\left\|u\right\|^{2}-{1\over{4d}}\left\|u\right\|^{2}_{L^{2}}-c\left\|u\right\|^{\alpha+1}_{L^{\alpha+1}}
\geq {1\over 4}\left\|u\right\|^{2}-c\eta^{\alpha+1}_{m-1}\left\|u\right\|^{\alpha+1}.$$
Let $r_{m}=[2(\alpha+1)c\eta^{\alpha+1}_{m-1}]^{{-1}\over{\alpha-1}}$, we have
$$f(u)\geq[{1\over 4}-{1\over{2(\alpha+1)}}]r^{2}_{m}=a_{m}\ with\ \left\|u\right\|=r_{m}.$$
Since $\alpha>1$, (3.12) shows that $a_{m}\longrightarrow\infty$ as $m\longrightarrow\infty$. $(f_{2})$ follows.\\
\end{Prf}
{\bf Lemma 3.3.} Let
$$g(u)=\int_{\rb}R(t,u)dt,\ \forall u\in E.$$
Then $g\in C^{1}(E,\rb)$ and $g'$ is a compact map.\\
\begin{Prf} By (3.3) and (3.4), $g\in C^{1}(E,\rb)$ and
$$g'(u)v=\int_{\rb}R'(t,u).vdt,\ \forall u,v\in E. \leqno(3.13)$$
Let $u_{n}\rightharpoonup u$ weakly in $E$. By Proposition 2.1, one can assume that $u_{n}\longrightarrow u$ strongly in $L^{p}$ for $p\in[1,\infty[$. By (3.13), we have
$$\left\|g'(u_{n})-g'(u)\right\|=\sup_{\left\|v\right\|=1}\left|\int_{\rb}(R'(t,u_{n})-R'(t,u)).v dt\right|.$$
By (3.3) there exists a constant $c_{1}>0$ such that for any $r>0$, one has
$$\left|\int_{\left|t\right|\geq r}(R'(t,u_{n})-R'(t,u)).vdt\right|
\leq c_{1}\int_{\left|t\right|\geq r}[\left|u_{n}\right|+\left|u\right|+\left|u_{n}\right|^{\alpha}+\left|u\right|^{\alpha}]\left|v\right|dt$$
$$\leq c_{1}[\left\|v\right\|_{L^{2}}(\int_{\left|t\right|\geq r}(\left|u_{n}\right|+\left|u\right|)^{2}dt)^{1\over 2} +\left\|v\right\|_{L^{\alpha+1}}(\int_{\left|t\right|\geq r}(\left|u_{n}\right|+\left|u\right|)^{\alpha+1}dt)^{\alpha\over{\alpha+1}}].$$
So by Proposition 2.1, there exists a constant $c_{2}>0$ such that for any $v\in E$, $\left\|v\right\|=1$,
$$\left|\int_{\left|t\right|\geq r}(R'(t,u(x_{n}))-R'(t,u)).vdt\right|$$
$$\leq c_{2}[(\int_{\left|t\right|\geq r}(\left|u_{n}\right|^{2}+\left|u\right|^{2})dt)^{1\over 2}+(\int_{\left|t\right|\geq r}(\left|u_{n}\right|^{\alpha+1}+\left|u\right|^{\alpha+1})dt)^{\alpha\over{\alpha+1}}]. \leqno(3.14)$$
We deduce from (3.14) that for any $\epsilon >0$, there exists $r>0$ so large such that
$$\left|\int_{\left|t\right|\geq r}(R'(t,u_{n})-R'(t,u)).vdt\right|<{\epsilon\over 2}, \leqno(3.15)$$
for all $n\in\nb$ and all $v\in E$, $\left\|v\right\|=1$. On the other hand, it is well known that since $u_{n}\longrightarrow u$, strongly in $L^{2}$,
$$\left\|R'(t,u_{n})-R'(t,u)\right\|_{L^{2}(I_{R})}\longrightarrow 0$$
as $n\longrightarrow\infty$, where $I_{r}=]-r,r[$. Therefore there is $n_{0}\in\nb$ such that
$$\left|\int_{\left|t\right|\leq r}(R'(t,u_{n})-R'(t,u)).v dt\right|<{\epsilon\over
2}, \leqno(3.16)$$
for all integer $n\geq n_{0}$ and all $v\in E$, $\left\|v\right\|=1$. Combining (3.15) and (3.16) yields
$$\left\|g'(u_{n})-g'(u)\right\|<\epsilon,\ \forall n\geq n_{0}.$$
Hence $g'$ is compact. \\
\end{Prf}
Finally, let us prove the Palais-Smale conditions.\\
{\bf Lemma 3.4.} $f$ satisfies $(PS)^{*}$ and $(PS)^{**}$ conditions.\\
\begin{Prf} The verification procedure for $(PS)^{*}$ and $(PS)^{**}$ conditions are the same, and so we only check the $(PS)^{*}$ condition. Suppose $u_{n}\in X_{n}$ be such that 
$$0\leq f(u_{n})\leq const.\ and\ f'_{n}(u_{n})\longrightarrow 0\ as\ n\longrightarrow\infty. \leqno(3.17)$$
We claim that $(u_{n})$ is bounded. If not, passing to a subsequence if necessary, we may assume that $\left\|u_{n}\right\|\longrightarrow\infty$ as $n\longrightarrow\infty$. By $(R_{4})$, $(R_{5})$, we have
$$2f(u_{n})-f'(u_{n}).u_{n}=\int_{\rb}[R'(t,u_{n}).u_{n}-2R(t,u_{n})]dt$$
$$\geq b\int_{\left\{t\in\rb/\left|u_{n}(t)\right|\geq r \right\}}\left|u_{n}\right|^{\beta}dt,\ \forall n\in\nb, \leqno (3.18)$$
which implies
$${1\over{\left\|u_{n}\right\|}}\int_{\left\{t\in\rb/\left|u_{n}(t)\right|\geq r\right\}}\left|u_{n}\right|^{\beta}dt\longrightarrow 0\ as\ n\longrightarrow\infty. \leqno(3.19)$$
By $(R_{2})$, $(R_{3})$ there exists a constant $c_{3}>0$ such that
$$\left|R'(t,x)\right|\leq c_{3}(\left|x\right|+\left|x\right|^{\alpha}),\ \forall x\in\rb^{2N}. \leqno(3.23)$$
Hence by (3.20)
$$f'_{n}(u_{n}).u^{+}_{n}={1\over 2}\left\|u^{+}_{n}\right\|^{2}-\int_{\rb}R'(t,u_{n}).u^{+}_{n}dt$$
$$\geq{1\over 2} \left\|u^{+}_{n}\right\|^{2} -\int_{\rb}\left|R'(t,u_{n})\right|\left|u^{+}_{n}\right|dt$$
$$\geq{1\over 2}\left\|u^{+}_{n}\right\|^{2} -c_{3}\int_{\rb}\left|u_{n}\right|^{\alpha}\left|u^{+}_{n}\right|dt-c_{3}\int_{\rb}\left|u_{n}\right|\left|u^{+}_{n}\right|dt. \leqno(3.21)$$
By H$\ddot{o}$lder's inequality and Proposition 2.1, we have
$$\int_{\rb}\left|u_{n}\right|^{\alpha}\left|u^{+}_{n}\right|dt= 
\int_{\left\{t\in\rb/\left|x_{n}(t)\right|\leq r\right\}}
\left|u_{n}\right|^{\alpha}\left|u^{+}_{n}\right|dt$$
$$+\int_{\left\{t\in\rb/\left|u_{n}(t)\right|\geq r\right\}}
\left|u_{n}\right|^{\alpha}\left|u^{+}_{n}\right|dt
\leq r^{\alpha} \int_{\left\{t\in\rb/\left|u_{n}(t)\right|\leq r\right\}}\left|u^{+}_{n}\right|dt$$       
$$+(\int_{\left\{t\in\rb/\left|x_{n}(t)\right|\geq r\right\}}
(\left|u_{n}\right|^{\alpha})^{\beta\over\alpha}dt)^{\alpha\over\beta}dt)
(\int_{\left\{t\in\rb/\left|u_{n}(t)\right|\geq r\right\}}
\left|u^{+}_{n}\right|^{\beta\over{\beta-\alpha}}dt)^{{\beta-\alpha}\over\beta}$$
$$\leq r^{\alpha} \left\|u^{+}_{n}\right\|_{L^{1}}+(\int_{\left\{t\in\rb/\left|u_{n}(t)\right|\geq r\right\}}\left|u_{n}\right|^{\beta}dt)^{\alpha\over\beta}\left\|u^{+}_{n}\right\|_{L^{\beta\over{\beta-\alpha}}}$$
$$\leq[\lambda_{1}r^{\alpha} +\lambda_{\beta\over{\beta-\alpha}}(\int_{\left\{t\in\rb/\left|u_{n}(t)\right|\geq r\right\}}\left|u_{n}\right|^{\beta}dt)^{\alpha\over\beta}]\left\|u^{+}_{n}\right\|.\leqno(3.22)$$
Similarly,
$$\int_{\rb}\left|u_{n}\right|\left|u^{+}_{n}\right|dt= 
\int_{\left\{t\in\rb/\left|u_{n}(t)\right|\leq r\right\}}
\left|u_{n}\right|\left|u^{+}_{n}\right|dt+
\int_{\left\{t\in\rb/\left|u_{n}(t)\right|\geq r\right\}}
\left|u_{n}\right|\left|u^{+}_{n}\right|dt$$
$$\leq r\left\|u^{+}_{n}\right\|_{L^{1}} +(\int_{\left\{t\in\rb/\left|u_{n}(t)\right|\geq r\right\}}\left|u_{n}\right|^{\beta}dt)^{1\over\beta}\left\|u^{+}_{n}\right\|_{L^{\beta'}}$$
$$\leq[\lambda_{1}r +\lambda_{\beta'}(\int_{\left\{t\in\rb/\left|u_{n}(t)\right|\geq r\right\}}\left|u_{n}\right|^{\beta}dt)^{1\over\beta}]\left\|u^{+}_{n}\right\|,\leqno(3.23)$$
where $\beta'$ is the H$\ddot{o}$lder's conjugate of $\beta$. Combining (3.21), (3.22), and (3.23), yields 
$$f'(u_{n}).u^{+}_{n}\geq{1\over 2} \left\|u^{+}_{n}\right\|^{2} -\lambda_{1}c_{3}[r^{\alpha}+ r]\left\|x^{+}_{n}\right\|$$
$$-c_{3}[\lambda_{\beta\over{\beta-\alpha}}(\int_{\left\{t\in\rb/\left|u_{n}(t)\right|\geq r \right\}}\left|u_{n}\right|^{\beta}dt)^{\alpha\over\beta}$$
$$+\lambda_{\beta'}(\int_{\left\{t\in\rb/\left|u_{n}(t)\right|\geq r \right\}}\left|u_{n}\right|^{\beta}dt)^{1\over\beta}]\left\|u^{+}_{n}\right\|. \leqno(3.24)$$
Since $1<\alpha<\beta$, we deduce from (3.19) and (3.24) that
$${{\left\|u^{+}_{n}\right\|}\over{\left\|u_{n}\right\|}}\longrightarrow 0\ as\ n\longrightarrow\infty. \leqno(3.25)$$
Similarly,
$${{\left\|u^{-}_{n}\right\|}\over{\left\|u_{n}\right\|}}\longrightarrow 0\ as\ n\longrightarrow\infty. \leqno(3.26)$$
Now, let
$$v_{n}(t)=\left\{
\begin{array}{l}
u_{n}(t),\ if\ \left|u_{n}(t)\right|\leq r,\\
0,\ if\ \left|u_{n}(t)\right|>r,
\end{array}\right.\leqno(3.27)$$
and
$$w_{n}(t)=u_{n}(t)-v_{n}(t)$$
for all integer $n$ and all $t\in\rb$. By (3.18) and (3.27), there exists a constant $c_{4}>0$ such that
$$c_{4}(1+\left\|u_{n}\right\|)\geq \left\|w_{n}\right\|^{\beta}_{L^{\beta}},\ \forall n\in\nb. \leqno(3.18)$$
Since $E^{0}$ is of finite dimension, we deduce from H$\ddot{o}$lder's inequality and (3.28)
$$\left\|u^{0}_{n}\right\|^{2}_{L^{2}}=<u^{0}_{n},u_{n}>_{L^{2}}=<u^{0}_{n},v_{n}>_{L^{2}}+<u^{0}_{n},w_{n}>_{L^{2}}$$
$$\leq r\left\|u^{0}_{n}\right\|_{L^{1}}+\left\|u^{0}_{n}\right\|_{L^{\beta'}}\left\|w_{n}\right\|_{L^{\beta}}\leq c_{5}\left\|u^{0}_{n}\right\|_{L^{2}}(1+\left\|w_{n}\right\|_{L^{\beta}}) \leqno(3.29)$$
for all integer $n$ and some constant $c_{5}>0$. Hence by (3.28) and (3.29), there exist positive constants $c_{6},\ c_{7}$ such that for all integer $n$
$$\left\|u^{0}_{n}\right\|\leq c_{6}\left\|u^{0}_{n}\right\|_{L^{2}}\leq c_{5}c_{6}(1+\left\|w_{n}\right\|_{L^{\beta}})\leq c_{7}(1+\left\|u_{n}\right\|^{1\over\beta}). \leqno(3.30)$$
Since $\beta>1$, we deduce from (3.30) that
$${{\left\|u^{0}_{n}\right\|}\over{\left\|u_{n}\right\|}}\longrightarrow 0\ as\ n\longrightarrow\infty. \leqno(3.31)$$
Combining (3.25), (3.26) and (3.31) yields
$$1={{\left\|u_{n}\right\|}\over{\left\|u_{n}\right\|}}\leq{{\left\|u^{-}_{n}\right\|+\left\|u^{0}_{n}\right\|+\left\|u^{+}_{n}\right\|}\over{\left\|u_{n}\right\|}}\longrightarrow 0\ as\ n\longrightarrow\infty, \leqno(3.32)$$
which is a contradiction. Hence $(u_{n})$ must be bounded. By Lemma 3.3, we deduce that $(u_{n})$ possesses a convergent subsequence, which completes the proof of Lemma 3.4. \\
\end{Prf}
The functional $f$ satisfies all the assumptions of Proposition 2.2, so it possesses a positive critical value sequence $(c_{k})$ satisfying $c_{k}\longrightarrow\infty$ as $k\longrightarrow\infty$. Therefore the system $(\mathcal H)$ possesses infinitely many homoclinic orbits $u_{k}$ satisfying $f(u_{k})\longrightarrow\infty$ as $k\longrightarrow\infty$.\\
{\bf References.} \par\noindent
[1] Y. Ding, "Existence and multiplicity results for homoclinic solutions to a class of Hamiltonian systems", Nonlinear Analysis, Vol. 25, No 11, pp 1095-1113, 1995. \par\noindent
[2] W. Omana and M. Willem, "Homoclinic orbits for a class of Hamiltonian systems", Diff. Integ. Equ., Vol.5, No 5, pp 1115-1120, 1992. \par\noindent
[3] S. Li, W. Zou, "Infinitely many homoclinic orbits for the second-order Hamiltonian systems", Applied Mathematics Letters 16, pp 1283-1287, 2003.\par\noindent
[4] X. Lin, X.H. Tang, "Infinitely many homoclinic orbits for Hamiltonian systems with indefinite sign subquadratic potentials", Nonlinear Analysis 74, pp 6314-6325, 2011.\par\noindent
[5] C. Liu, Q. Zhang, "Infinitely many homoclinic solutions for second order Hamiltonian systems", Nonlinear Analysis 72, pp 894-903, 2010.\par\noindent
[6] Y. Lv, C.L Tang, "Existence of even homoclinic orbits for second-order Hamiltonian systems" Nonlinear Analysis 67, pp 2189-2198, 2007.\par\noindent
[7] Z.Q. Ou, C.L. Tang, "Existence of homoclinic solution for the second order Hamiltonian systems", J. Math. Anal. Appl. 191, pp 2003-2011, 2004.\par\noindent
[8] P.Rabinowitz, "Homoclinic orbits for a class of Hamiltonian systems", Proc.Roy.Soc. Edinburgh Sect. A 114, pp 33.38, 1990. \par\noindent
[9] C.L. Tang, L.L. Wan, "Existence of homoclinic orbits for second order Hamiltonian systems without (AR) condition", Nonlinear Analysis 74, pp 5303-5313, 2011.\par\noindent
[10] J. Yang, F. Zhang, "Infinitely many homoclinic orbits for the second order Hamiltonian systems with super-quadratic potentials", Nonlinear Analysis: Real World Applications 10, pp 1417-1423, 2009.\par\noindent
[11] Y. Ding, "Infinitely many homoclinic orbits for a class of Hamiltonian systems with symmetry", Chin. Ann. of Math. 19B: 2, pp 167-178, 1998.\par\noindent
[12] Y. Ding, M. Girardi, "infinitely many homoclinic orbits of a Hamiltonian system with symmetry", Nonlinear Analysis 38, pp 391-415, 1999.\par\noindent
[13] Y. Ding, L. Jeanjean, "Homoclinic orbits for a nonperiodic Hamiltonian system", J. Differential Equations 237, pp 473-490, 2007.\par\noindent
[14] Y. Ding and S. Li, "Homoclinic orbits for first order Hamiltonian systems", J. Math. Anal. Appl. 189, pp 585-601, 1995. \par\noindent
[15] I. Ekeland, V. Coti-Zelati, "A variational approach to homoclinic orbits in Hamiltonian systems", Math. Ann. 288, pp 133-160, 1990. \par\noindent
[16] K. Tanaka, "Homoclinic orbits in a first order superquadratic Hamiltonian system: covergence of subharmonic orbits", J. Diff. Equ. 94, pp 315-339, 1991. \par\noindent
[17] M. Timoumi, "On homoclinic orbits for a class of nocoercive superquadratic Hamiltonian systems", Nonlinear Analysis 74, pp 5892-5901, 2011\par\noindent
[18] D.E. Edmunds and W.D. Evans, "Spectral theory and differential operators", Clarendon, Oxford, 1987. \par\noindent
[19] G. Fei, "On periodic solutions of superquadratic Hamiltonian systems", Electron. J. Differential Equations 8, 12 p, 2002.\par\noindent
[20] M. Timoumi, "Periodic solutions for noncoercive superquadratic Hamiltonian systems", Demonstratio Mathematica Vol. XL, No 2, pp 331-346, 2007.\par\noindent

\end{document}